\documentclass[12pt]{article}
\usepackage{fancyhdr}

\usepackage{amsmath}
\usepackage{mathtools}
\usepackage{amsfonts}
\usepackage{amsthm}
\usepackage{amssymb}
\usepackage{extarrows}
\usepackage{color}
\usepackage{dsfont}
\usepackage{geometry} 
\newgeometry{tmargin=2.7cm, bmargin=3.5cm, lmargin=2.7cm, rmargin=2.8cm} 

\usepackage{color}

\title{A pocket guide to nonlinear differential equations\\ in Musielak--Orlicz spaces}

\usepackage{authblk}

\author[$\dagger$]{Iwona Chlebicka (Skrzypczak)\thanks{email address: \texttt{ iskrzypczak@mimuw.edu.pl}\\ The author is supported by NCN grant no. 2016/23/D/ST1/01072.
}}
\affil[$\dagger$]{\small
Institute of  Mathematics, Polish Academy of Sciences, \newline
ul. \'{S}niadeckich 8, 00-656 Warsaw, Poland
}

\newcommand{\barint}{
         \rule[.036in]{.12in}{.009in}\kern-.16in
          \displaystyle\int  }
\def\R{{\mathbb{R}}}

\def\r{{\mathbb{R}}}
\def\N{{\mathbb{N}}}
\def\rn{{\mathbb{R}^{N}}}

\def\w{\widetilde}

\def\rp0{{[0,\infty)}}

\def\ve{{\varepsilon}}

\newtheorem{theo}{\bf Theorem}   
\newtheorem{rem}{\bf Remark} 
\newtheorem{defi}{\bf Definition} 
\newtheorem{ex}{\bf Example} 
\newtheorem{fact}{\bf Fact} 
\newtheorem{prop}{\bf Proposition} 
\newtheorem*{op}{\bf Open problem}

\newcommand{\vp}{\varphi}
\newcommand{\pa}{\partial}
\newcommand{\dv}{\mathrm{div}}

\date{}
\begin{document}
\maketitle \sloppy

\thispagestyle{empty}


\parindent 1em

\begin{abstract}

The Musielak--Orlicz setting unifies variable exponent, Orlicz, weighted Sobolev, and double-phase spaces. They inherit technical difficulties resulting from general growth and inhomogeneity.

In this survey we present an overview of developments of the theory of existence of~PDEs in the setting including reflexive and non-reflexive cases, as well as isotropic and anisotropic ones. Particular attention is paid to problems with data below natural duality in absence of Lavrentiev's phenomenon.

\end{abstract}

\bigskip

  
  \medskip


\newpage

\section{Introduction}

Vast literature describes various aspects of PDEs with the leading part of the operator   having a~power--type growth with the preeminent example of the $p$-Laplacian. There is a~wide range of directions in~which the polynomial growth case has been developed including variable exponent, convex, weighted, and double-phase approaches. Musielak-Orlicz spaces  cover all of the mentioned cases. They are equipped with the norm defined by~means of 
\[\varrho_M(\nabla u)=\int_\Omega M(x,\nabla u)\,dx,\]
called \textbf{modular} and $M:\Omega\times\rn\to\rp0$ is then called a {modular function}. We call this function non-homogeneous due to~its $x$-dependence. We assume $M$ to be convex with respect to the gradient variable. Note that $M$ can be an anisotropic function, which means that it depends on whole gradient $\nabla u$, not only its length $|\nabla u|$ and the behaviour of~$M$ may vary in different directions. The typical anisotropic example is~$M(x,\nabla u)=\sum_{i=1}^N |u_{x_i}|^{p_i(x)}$, but the anisotropy does not have to be expressed by separation of roles of~coordinates. 

Let us point out a few basic references. We refer to pioneering works  by Orlicz~\cite{Orlicz31}, where variable exponent spaces are introduced, and by Zygmund~\cite{Zygmund} with elements of~Orlicz spaces. Variable exponent spaces are carefully treated in the monographs by Cruz-Uribe and Fiorenza~\cite{cruz-u} and Diening, Harjulehto, H\"{a}st\"{o}, and R{{u}}{\v{z}}i{\v{c}}ka~\cite{ks}. The foundations of Orlicz spaces are described by Krasnosel'skii and Rutickii~\cite{KraRu} and Rao and Ren in~\cite{rao-ren}. The most exhaustive study on weighted Sobolev setting is presented by Turesson in~\cite{Turesson}. 
 
Passing to general Musielak-Orlicz spaces -- the first systematic approach to non-homogeneous setting with general growth was provided by Nakano~\cite{Nakano}, then Skaff~\cite{Sk1,Sk2} and Hudzik~\cite{HH3,HH4}, but the key role in the functional analysis of  modular spaces is played by the comprehensive book by~Musielak~\cite{Musielak}. None of these sources is focused on the theory of PDEs. 
  
  \medskip
  
 We start with a brief presentation  of spaces included in the framework of~the~Musielak-Orlicz setting in connection to PDEs highlighting only samples of vastness of~results therein. Our aim is to present difficulties that each of examples carry rather than to~provide a comprehensive overview of results in each setting.   Moreover, to~keep our guide pocket-sized we restrict ourselves to the  theory of existence of solutions to~problems in classical euclidean spaces slipping over regularity issues comprehensively described in~\cite{min-dark} with further developments in variable exponent spaces mentioned in~\cite{overview}. Particular attention is paid here to~the issue of existence to~problems with data below natural duality and relate them either to growth conditions, or to~the~absence of~Lavrentiev's phenomenon when asymptotic behaviour of a modular function is sufficiently balanced. 

\section{Overview of spaces}
 This section is devoted to concise summary of features and difficulties of spaces included in the framework of the Musielak-Orlicz setting.

\subsection{Sobolev  spaces: classical, weighted, and anisotropic}
The natural setting to study solutions to elliptic and parabolic partial differential equations involving Laplace or $p$-Laplace operator
\[\Delta_p u=\dv \big(|\nabla u|^{p-2}\nabla u\big)\] are classical Sobolev spaces, when the modular function has a form $M(x,\nabla u)=|\nabla u|^{p}$. Let us restrict ourselves to the classical references~\cite{ladyzenska,ladyzenska-ell}.  However, if one wants to~study more degenerate partial differential problems involving various types of singularities in the coefficients, e.g. the weighted $\omega$-$p$-Laplacian
\[\Delta_p^\omega u=\dv \big(\omega(x)|\nabla u|^{p-2}\nabla u\big),\]
 the relevant setting are weighted Sobolev spaces,  see~\cite{FaKeSe,HeKiMa}, where the modular function is~$M(x,\nabla u)=\omega(x)|\nabla u|^{p}$.

 We refer to~\cite{kuf-opic}, where Kufner and Opic introduce the assumption sufficient for the~weighted space to be continuously embedded in $L^1_{loc}(\Omega),$ and consequently for any function from the weighted space to have well-defined distributional derivatives. This condition is called $B_p$--condition, and yields that the weight $\omega$ is a~positive a.e.  Borel  measurable function such that $
 \omega'=\omega^{-1/(p-1)}(x)\in L^1_{{  loc}}(\Omega)$. This condition is weaker than $A_p$-condition, cf.~\cite{muckenhoupt}.  
Turesson's book~\cite{Turesson} consists of a comprehensive study on~the case of $A_p$--weights. It provides weighted analogues of~multiple results from the theory of~non-weighted Sobolev spaces and from non-weighted potential theory. Useful embeddings are proven in~\cite{anh,OpKuf}. PDEs in the weighted setting are considered e.g. in~\cite{weightmatteo1,nikita,CalMu,pdakis1,isazg2,FaJeKe,FaKeJe,
FeGrMe}.

\bigskip

To describe anisotropy, that is when energy density is not the same in various directions, one can use the Sobolev space with different exponents. The model example involves anisotropic $\vec{p}$-Laplacian
\[\Delta_{\vec{p}} u=\dv\left(\sum_{i=1}^N |u_{x_i}|^{p_i-2}u_{x_i}\right)\quad\text{with}\quad\vec{p}=(p_1,\dots,p_N)\]
and thus, the appropriate space is
\[
W^{1,\vec{p}}(\Omega):= \Big\{ f\in W^{1,1}_{loc}(\Omega) : \ \ \  f\in L^{p_0}(\Omega),\ \  
 f_{x_i}\in  L^{p_i} (\Omega),\ \ \text{for }\ i=1,\dots,N\Big\},
\]
where $p_0$ is a harmonic mean of $p_1,\dots,p_N$. See~\cite{Ci-fully} for the embedding result.

The research on partial differential equations in the {anisotropic} setting comes back to~\cite{Ma-sym,Ta0,UrUr,We-sym} and with merely integrable data to~\cite{Aniso-BGM}. For some results on~regularity and existence we refer to e.g.~\cite{ELM,FraGaKa,FraGaLi,Lieb-aniso,BStr,Ve1}, for other estimates on~solutions to~\cite{AlBlFe,AlCi-aniso,AlBlFe-p}, while for nonexistence to~\cite{ambr09,FraGaKa}.

\subsection{Variable exponent spaces}\label{ssec:varexp}

To describe setting in which energy density varies with the space variable, the basic idea is to consider
\[\Delta_{p(x)} u=\dv \Big(|\nabla u|^{p(x)-2}\nabla u\Big)\quad\text{or}
\quad \w{\Delta}_{p(x)} u=\dv \Big(p(x)|\nabla u|^{p(x)-2}\nabla u\Big).\]
Therefore, the relevant setting is provided in the framework of the variable exponent Lebesgue space  defined by the mean of the modular function $M(x,\nabla u)=|\nabla u|^{p(x)}$. Namely,
\[
    L^{p(\cdot)}(\Omega)=\left\{f \text{ - measurable } : \ \ \  \int_{\Omega} |f(x)|^{p(x)} dx<\infty \right\}.
\]

Variable exponent Lebesgue and Sobolev spaces have received significant attention. For good theoretical basis we refer to books of Cruz-Uribe and Fiorenza~\cite{cruz-u} and of~Diening, Harjulehto, H\"{a}st\"{o}, and R{{u}}{\v{z}}i{\v{c}}ka~\cite{ks}. Typical examples of  equations involving the variable exponent setting include models of electrorheological fluids~\cite{mingione02,raj-ru,el-rh2}, image restoration processing~\cite{chen06}, elasticity equations~\cite{zhikov97}, and thermistor model~\cite{zhikov9798}. Since the setting is carefully examined and broad range of problems is already addressed, we mention only a~few attempts to basic properties of PDEs such as existence~\cite{din06,fanzhang,liu,puczha14}, regularity results~\cite{mingione06,mingione01}, maximal principle~\cite{HHLT}, and nonexistence~\cite{adagor14,sdis}. Existence to~$L^1$-data problems within isotropic approach were studied in~\cite{BWZ,wit-zim} and  anisotropic case in~\cite{BeKaSa,BeLaSa}. For more information  on this issue see Section~\ref{sec:below}. Let us refer also to~the~survey~\cite{overview}  describing carefully developments of the theory of~differential equations within the setting.
Complicated further setting of weighted variable exponent spaces are considered e.g. from the point of~view of harmonic analysis in~\cite{CUDH}.

\medskip

Variable exponent spaces are reflexive for every exponent $1<<p<<\infty$, which equippes with lots of tools of functional analysis. Nonetheless, we deal with difficulties resulting from non-homogeneity of such space, namely, from the fact that density of~smooth functions depends on the regularity of a modular function. If exponent is not regular enough, we meet so-called Lavrentiev's phenomenon~\cite{LM,ZV}, originally meaning the situation when the infimum of a~variational problem over regular functions (e.g. smooth or Lipschitz) is strictly greater than infimum taken over the set of all functions satisfying the~same boundary conditions. See e.g.~\cite[Example~3.2]{ZV} by Zhikov with $p(\cdot)$ being a step function. The notion of the Lavrentiev phenomenon became naturally generalised to~describe the situation, where functions from certain spaces cannot be approximated by regular ones, cf. Section~\ref{sec:appr}.

Its presence plays an important role in the calculus of variations~\cite{min-double-reg1}, but also in the~existence theory~\cite{ks}, and homogenization~(see Zhikov's pioneering work~\cite{ZV} and more in~\cite{Zhikov2011,ZKO}). Typically to ensure density of smooth functions the assumption imposed on exponent is log-H\"older continuity, see~\cite{ks} where one can find an explanation that this is not always a necessary condition. Nonetheless, it is handy and certainly the most commonly assumed one. This entails that, although the $x$-dependence of $p$ enables to study problems with different exponents in various subdomains, PDEs considered in the setting are in the overwhelming majority formulated for log-H\"older exponents, which excludes dramatic changes of the  energy density. The opposite approach is taken into account in double-phase spaces.

\subsection{Double-phase spaces}\label{ssec:double-phase}

The natural direction of relaxing power-type growth of the leading part of the operator relaxing  is considering  the growth sandwiched between two power-type functions  comes from Marcellini~\cite{Marc1,Marc2}. Here we concentrate on the particular non-uniformly elliptic operator of this kind, namely, \[\dv\Big(\big(|\nabla u|^{p-2}+a(x)|\nabla u|^{q-2}\big)\nabla u\Big)\]  with $1 < p<q < \infty $ and a weight function $a:\Omega\to[0,\infty)$ which may disappear, that can describe the diffusion-type process in the space whose certain subdomains are distinguished from the others. In this case we shall consider modular function given by the formula\[ M(x,\nabla u)=H(x,|\nabla u|)=|\nabla u|^p+a(x)|\nabla u|^q 
\] 
and -- corresponding to solutions to an equation involving the above non-uniformly elliptic operator -- minimizers of variational functionals, cf.~\cite{bcm17,min-double-reg2,min-double-reg1,colombo}.

 This case is related to the variable exponent space with an exponent being a~step function, rather than with the weighted Sobolev space. With this analogy notice that it describes the composite material having on $\{x\in\Omega:\,a(x)=0\}$ the energy density with $p$-growth, whereas on $\{x\in\Omega:\,a(x)>0\}$ the growth of order $q$. The relation between the double-phase space and the variable exponent one is explored in~\cite{bcm-st}. The space investigated therein, being kind of borderline line case of the double-phase space, is equipped with\[ M(x,\nabla u)=H_{\log}(x,|\nabla u|)=|\nabla u|^p+a(x)|\nabla u|^p\log(e+|\nabla u|).\]
The key feature of both settings (the double-phase space and its borderline case) is that regularity of the weight function $a$ dictates the ellipticity
rate of~the energy density indicating the range of parameters necessary to ensure good properties of~the space including the modular approximation, cf. Definition~\ref{def:convmod} and Remark~\ref{rem:double-phase}.  

The double-phase spaces appeared originally in context of homogenization and the~Lavrentiev phenomenon (see~\cite{Zhikov2011} and Section~\ref{ssec:varexp}). Let us mention that recently regularity theory of minimisers to variational functionals in this setting is getting increasing attention, starting from~\cite{min-double-reg1,min-double-reg2}. In this context the optimal approximation in the modular topology (called also absence of Lavrentiev's phenomenon in~\cite{min-double-reg1} and Gossez's approximation theorems in~\cite{yags}, cf. Section~\ref{sec:appr}) plays crucial role in the sharp regularity results. As far as a variable exponent is expected to be log-H\"older continuous, here the exponents $p<q$ should be close to each other. The optimal closeness condition in the case of a~H\"older continuous weight $a\in C^{0,\alpha}$ obtained in~\cite{ELM} requires $q/p\leq 1+\alpha/N$. There is an~interesting counterexample provided by Colombo and Mingione in~\cite[Theorem~4.1]{min-double-reg1} that illustrates the sharpness of this condition. To exclude Lavrentiev's phenomenon in the borderline case it is needed that $a$ is log-H\"older continuous -- the same what is assumed on the variable exponent in the corresponding context. What is more, it is interesting to note that, while considering the variable exponent functional and the borderline case of the double phase functional, the regularity conditions to assume on the exponent $p(x)$ and the coefficient $a(x)$ in order to get regularity of minimizers, are the same. See~\cite{bcm-st} for detailed information on the nature of~the mutual relation of these two function spaces.
 
 Let us note that the double-phase space together with the borderline case with bounded $a\geq 0$ and $1 < p,q < \infty $ are reflexive, no matter how is the interplay of the parameters or~how irregular the weight is.

\subsection{Isotropic and anisotropic Orlicz spaces}\label{sec:orlicz}

Let us concentrate on the operator of the form \[\dv\left(\frac{B(\nabla u)}{|\nabla u|^2}\nabla u\right),\]
where to come back to $p$--Laplacian one shall choose $B(t)=|t|^p$. Problems involving $B$ with growth more general  than just of power-type are studied from Talenti~\cite{Ta1}, Donaldson~\cite{Donaldson}, and Gossez~\cite{Gossez} with new spirit given by Fusco and Sbordone~\cite{Fusco-Sbordone}, Lieberman~\cite{Lieb91} and Cianchi~\cite{Ci97}. The applied motivation for the Orlicz setting includes modelling of non-Newtonian fluids~\cite{BrCi} and of~elastodynamics~\cite{NgPa}. In the framework of Orlicz spaces a modular function then has the form $M(x,\nabla u)=B(\nabla u)$ with $B$ being an $N$-function, i.e. convex function with $B(0)=0$ and accelerating faster than linear in the origin and in the infinity. The natural setting is the space
\[
    L_{B}(\Omega)=\left\{f \text{ - measurable } :\ \ \ \int_{\Omega}B( f(x)) dx<\infty \right\}.
\]
 
As in the variable exponent case the relevant form of norm has the Luxemburg type, i.e. 
\[    \|f\|_{L_{B}(\Omega)}:=\inf \Big\{\lambda>0: \int_{\Omega} B\left(\frac{f(x)}{\lambda}\right) dx \leq 1 \Big\}.\]
Note that this definition is not isotropic. This approach is classically homogeneous ($B$ is not a function of the space variable), but admit modular functions $B$ far from being of~power-type as well. In particular, we can expect 
exponentially fast growth \[B(t)=|t|\exp| t|\] or slower than any power strictly bigger than $1$ \[B(t)=|t|\log(1+|t|)\] or mixed condition (other in the origin and in infinity, or in various directions).

\medskip

We refer to the short book with the basics of isotropic Orlicz spaces by Krasnosel'skii and Rutickii~\cite{KraRu} introducing the setting in a geometric way. The classical,  comprehensive book of Rao and Ren~\cite{rao-ren} systematises the framework, while~\cite[Section~8]{adams-fournier} highlights clearly the crucial points of the theory relevant to isotropic differential equations. 

\medskip

In the functional analysis of the Orlicz setting the key role is played by  $B^\ast$ -- the defined below  complementary~function (called also the Young conjugate, or the Legendre transform) to a function  $B:\rn\to\r$. The complementary~function is given by the following anisotropic formula
\[B^\ast(\eta)=\sup_{\xi\in\rn}(\xi\cdot\eta-B(\xi)),\qquad \eta\in\rn.\]
See~\cite{KraRu} for appropriate explanation how $B$ and $B^\ast$ complement each other accompanied by nice graphs. To give a basic example let us note that
\[B(s)=\frac{1}{p}|s|^p\qquad\implies\qquad B^\ast(s)=\frac{1}{p'}|s|^{p'}.\]
Typically the assumptions on growth control of $B$ play technical role, whereas describing growth of $B^\ast$ fixes a structure of a space.
 
\subsubsection*{Growth restrictions}
Major part of the studies concern the case when a modular function satisfies $\Delta_2$--condition called also the doubling condition and denoted $M\in\Delta_2$ or $\Delta_2(\{M\})<\infty$. We say that an $N$-function $M:\rn\to\r$ satisfies $\Delta_2$--condition, if~there exists a constant $c_{\Delta_2}>0$ such that
\[ B(2\xi)\leq c_{\Delta_2}B(\xi).
\] 
Note that the condition is anisotropic, since it describes the behaviour in every direction.

Obviously, power-type functions \[B_1(\xi)=|\xi|^p\quad \text{with}\quad 1<p<\infty,\] as well as Zygmund-type functions \[B_2(\xi)=|\xi|^p \log^\alpha (1+|\xi|) \quad \text{with}\quad 1<p<\infty\ \text{and}\ \alpha\geq 0\] satisfy   $\Delta_2$-condition. Moreover, their conjugates $B_1^\ast,B_2^\ast\in\Delta_2$. What is more, the family of~functions satisfying $\Delta_2$-condition is invariant with respect to multiplications and compositions. Let us indicate that this condition excludes too quick growth. When we consider \[B(\xi) = (1+|\xi|)\log(1+|\xi|)-|\xi|\quad\in\quad\Delta_2,\] we note that its complementary function does not satisfy $\Delta_2$-condition. Indeed, \[B^\ast(\xi ) = \exp(|\xi|)-|\xi|-1 \quad\not\in\quad\Delta_2. \]
 
Moreover, in the isotropic case when a function $B$ is differentiable the condition \[ 1<i_B=\inf_{t>0}\frac{tB'(t)}{B(t)}\leq \sup_{t>0}\frac{tB'(t)}{B(t)}=  s_B<\infty
\] is equivalent to $B, {B}^\ast\in\Delta_2$,~\cite[Section~2.3, Theorem~3]{rao-ren}, denoted sometimes $\Delta_2(\{B,B^\ast\})<\infty$, cf.~\cite{DieEtt}. Indeed, if $s_B<\infty$ then $B\in\Delta_2$, whereas $i_B>1$ entails the $\Delta_2$-condition imposed on $ {B}^\ast$.  Note that if $B,B^\ast\in\Delta_2,$ then $B$ is trapped between two power--type functions, i.e.
\[\frac{B(t)}{t^{i_B}}\quad\text{is nondecreasing }\qquad\text{and}\qquad \frac{B(t)}{t^{s_B}}\quad\text{is nonincreasing,}
\]
 see~\cite[Section~2.3]{rao-ren}. The reverse implication is not true. There exist  functions sandwiched between $t^p$ and $t^q$ for arbitrary $1<p<q$ and not satisfying $\Delta_2$-condition. See~\cite{CGZG} for a~construction.  

\medskip

If we do not restrict the growth of a modular function $B$, the Orlicz space $L^B$ is not reflexive and -- thus -- weak and weak-$\ast$ convergence does not coincide. What is more, in reflexive spaces we are equipped with Mazur's Lemma implying that we do not leave the space, when we construct a function as a weak limit of smooth ones. To ensure reflexiveness in the Orlicz case, it is necessary and sufficient to impose $\Delta_2$-condition on~both -- a modular function and its conjugate. See related Section~\ref{ssec:genMOs}. 

\subsubsection*{Modular topology}

 When the growth of~a modular function is arbitrary another type of approximation is more appropriate. In his seminal paper~\cite{Gossez}  Gossez proves that weak derivatives in Orlicz-Sobolev spaces are strong derivatives with respect to the modular topology given by the definition of~the modular convergence   
\[\xi_\delta\xrightarrow[\delta\to 0]{M}\xi\qquad\iff\qquad\exists_{\lambda>0}:\qquad
\int_{\Omega}B\left( \frac{\xi_\delta-\xi}{\lambda }\right)dx\to 0.\]
Therefore, the modular closure is the relevant tool in this setting.

\subsubsection*{Embeddings}

For Sobolev--Orlicz spaces expected embedding theorems hold true, namely
\[W_0^{1,B}(\Omega)\hookrightarrow{} L_{\hat{B}}(\Omega),\qquad \Omega\subset\rn,\]
with $\hat{B}$ growing in a certain sense faster than $B$.  The classical references are~\cite{Ad-emb,DT-emb,Ta89,Trud}. The embedding theorems optimal from certain points of view comes from Cianchi and are concisely described in~\cite{Ci-aspects}. In the case of quickly growing modular function, i.e. under the integral condition corresponding to the case of $p$-growth with $p>N$, it holds that $W_0^{1,B}\hookrightarrow{}L^\infty$. When $B$ grows slowly, corresponding to $p\leq N$, we expect $W_0^{1,B}\hookrightarrow{} L_{\hat{B}}.$ Let us refer to~\cite{Ci96-emb} for best possible embeddings into Orlicz target space and to~\cite{cianchi04} for embeddings into Lorentz-Orlicz spaces being optimal target among all rearragement invariant spaces. The fully anisotropic optimal version of~\cite{Ci96-emb} is provided in~\cite{Ci-fully}, whereas the higher-order embeddings can be found in~\cite{Ci-high-or}.
 
There are certain simple versions of modular Sobolev-Poincar\'e inequalities implying $W^{1,B}_0 (\Omega)\hookrightarrow{} L_{B^{N'}}(\Omega)$, $N'=N/(N-1)$. They are indeed far from optimal, but for some purposes very handy since capturing all kinds of growth, see~\cite{Baroni-Riesz,CGZG}. 

\subsubsection*{Sample of PDE results} The classical reference for existence in the reflexive isotropic Orlicz-Sobolev setting is the already mentioned paper~\cite{Ta1} by Talenti. We refer for other results on existence and gradient estimates to~\cite{ACCZG,AlCiSb,CiMa,Cap,IC-grad,CiMa-ARMA14,CiMa-ARMA18}, nonexistence~\cite{FiGi,nonex}, regularity~\cite{Lieb91,Lieb93,CiFu} and then \cite{BrVe,DieEtt}, and foundations of the potential theory to~\cite{Baroni-Riesz}.

The cornerstones of nonlinear boundary value problems in~non-reflexive Orlicz-Sobolev-type setting were laid by Donaldson~\cite{Donaldson} and  Gossez~\cite{Gossez2,Gossez3,Gossez}. We refer to a very nice survey on elliptic problems~\cite{Mustonen} by Mustonen and Tienari, while the most relevant reference on parabolic ones are~\cite{ElMes,ElMes2} by Elmahi and Meskine. This research was continued in~direction of problems with data below natural duality~\cite{ElMes-ellL1,ElMes-ellL12} and the degree theory~\cite{GoMu,LaMu,Tie}. Fully anisotropic Orlicz spaces were engaged in PDEs in~\cite{Ci-sym} and~\cite{Alb-CPDE11,ACCZG,AlCi-aniso,AlBlFe-fully,BarCi-aniso}.

Amongst the very few regularity results outside spaces when a modular function is controlled by power-type functions, we shall refer to~\cite{Fuchs-Mingione} for investigations on slowly growing functionals ($L\log L$), to \cite{Marc-everywhere} capturing quick growth (exponential), and~\cite{marc-papi} admitting both slow and quick (but not arbitrary) growth.  We give more information on solutions to problems with data below duality in Section~\ref{sec:below}, thus we give here sole references to~\cite{renel2,renel3,le-ex}.

\subsubsection*{Difficulties} 

Let us describe main struggles of the Orlicz setting with some ideas how to cope with them.

\medskip

\textit{\underline{Anisotropy}. } Anisotropic modular function that can be considered is \textbf{not} necessarily  describing each direction separately, i.e.
\[B(\xi)=\sum_{i=1}^NB_i({\xi_i})\qquad\xi=(\xi_1,\dots,\xi_N)\in\rn.\]
A modular function (and then the space) which does not admit the above decomposition is called \textbf{fully anisotropic}. The two-dimensional example of fully anisotropic function provided in~\cite{Trud-Ex} is\[B(\xi)=|\xi_1-\xi_2|^\alpha+|\xi_1|^\beta\log^\delta(c+|\xi_1|),\qquad \alpha,\beta\geq 1,\] and $\delta\in \R$ if $\beta>1$, or $\delta>0$ if $\beta=1$, with $c>>1$ large enough to ensure convexity. 

In the fully anisotropic case to ensure that the Luxemburg-type norm is still a norm, besides convexity, a modular function has to be even ($B(\xi)=B(-\xi)$). See Definition~\ref{def:Nf} and ignore $x$-dependence. Note that the structure of~a space is poorly controlled and admissible tools are restricted. The behaviour of $B$ can be wilder than in the above example, because the speed of growth can change dramatically with the direction. Let us point out that unlike the isotropic case, where we can define the Orlicz-Sobolev space as
\[\left\{u\in L_B(\Omega):\ \nabla u\in L_B(\Omega;\rn)\right\},\]
 in the fully anisotropic case, we shall rather consider 
\[\left\{u\in L^1(\Omega):\ \nabla u\in L_B(\Omega;\rn)\right\}\qquad\text{or}\qquad \left\{u\in L_b(\Omega):\ \nabla u\in L_B(\Omega;\rn)\right\}\]
with some isotropic $N$-function $b$ growing essentially less rapidly than $B$. Note that symmetrization technique~\cite{Ci-fully} leads to a different embedding i.e. $u\in L_b$ with $b=B_N$ constructed therein.

Many facts holding in the isotropic case are not true in the anisotropic setting anymore. It concerns in particular those which describe the interplay between a behaviour of a~modular function and its complementary. 

The technique with the fundamental meaning in anisotropic Orlicz spaces is symmetrization coming from~\cite{Kl-sym0,PoSze}. It has been developed in PDE context in~\cite{Ci-fully,Ci-sym} and applied in studies on existence and regularity~\cite{Alb-CPDE11,AlBlFe,AlBlFe-fully,ACCZG,AlCi-aniso,BarCi-aniso}. For symmetrization-free approach in the setting let us refer to~\cite{Gparabolic,le-ex} and results in the Musielak-Orlicz setting mentioned in Sections~\ref{ssec:genMOs} and~\ref{sec:below}. The important issue is then to choose proper growth and coercivity conditions to by-pass the lack of~structure of a space.

 \medskip

\textit{\underline{Description of growth and coercivity}. } Let us consider problems modelled upon \[-\dv A=f\qquad\text{ or }\qquad u_t-\dv A=F,\] where $A$ is going to be governed by a modular function $B$.

In an Orlicz space \[B,B^\ast\in\Delta_2\qquad\iff\qquad L_{B^\ast}=(L_B)^\ast,\] see~\cite[Section~8]{adams-fournier}. See also Section~\ref{ssec:genMOs} for Definition~\ref{def:MOsp} and some remarks on the functional analysis of Musielak-Orlicz spaces, that applies also here while ignoring $x$-dependence. Note that in the Orlicz case classes $E_B$ and $L_{B}$ does not coincide outside $\Delta_2$--family of~the modular functions.

\medskip
 
 Thus, growth and coercivity conditions can be expressed similarly to the typical assumptions in the power-type case with $B(\xi)=|\xi|^p$, $1<p<\infty$, namely
 \[A(x,\xi)\xi\geq B(|\xi|)\qquad\text{and}\qquad |A(x,\xi)|\leq c B'(|\xi|),\]
cf.~\cite{CiMa}, which is enough to control the structure of the space.  The reason for which the above conditions are sufficient is the fact that in the doubling case we have \[B^\ast\big(B'(s)\big)\leq c B(s)\qquad\text{implying}\qquad B'(s) \leq (B^\ast)^{-1}\big(c B(s)\big).\]
Then we deal with proper dual pairing. Namely, for a solution $u$, such that $\nabla u\in L_B$ we get $A(\cdot,\nabla u)\in L_{B^\ast}=(L_B)^\ast$. Let us stress again that it is not the case of $B$ without growth restrictions, e.g. spaces $L\log L$ or $L_{\exp}$.   In fact, outside the doubling case we shall rather consider \[A(x,\xi)\xi\geq B(|\xi|)\qquad\text{and}\qquad |A(x,\xi)|\leq c_1(B^\ast)^{-1}\big(c_2 B(|\xi|)\big),\] 
cf.~\cite{ElMes,Gossez2,Cap,Mustonen}.

If one wants to keep anisotropic structure  more relevant would be to~consider
\[A(x,\xi)\xi\geq B( \xi )\qquad\text{and}\qquad B^\ast\big(c_1A(x,\xi)\big)\leq c_2 B( \xi )\]
or to hide them both in one assumption as e.g. in~\cite{gwiazda-ren-ell,pgisazg1}
\[c\Big( B( \xi)+B^\ast\big(A(x,\xi)\big)\Big)\leq A(x,\xi)\xi.\]

\medskip

\textit{\underline{Lack of factorization}. } On the other hand, there are noticable difficulties in the study on parabolic problems resulting from the lack of the integration by parts formula, which therefore has to be formulated in an advanced way. There are two main reasons  for this. The first one is the lack of strong density of smooth functions described more carefully in Section~\ref{sec:appr} (as already mentioned above the modular topology in more relevant in the setting), whereas the second one is that in general\[L_B(\Omega_T)\neq L_B(0,T;L_B(\Omega)).
\]
Indeed, to ensure equality, the growth of $B$ has to be restricted.
\begin{prop}[\cite{ZaGa}] Suppose $\Omega\subset\rn$ is a bounded subset, $T<\infty$, $B:\rp0\to\rp0$ is an $N$-function. Then
\[ L_B(\Omega_T)= L_B(0,T;L_B(\Omega))\]
if and only if there exist constants $k_0,k_1$, such that for every $s\geq 1/T$ and $r\geq 1/|\Omega|$ it holds that
\[k_0B^{-1}(s)B^{-1}(r)\leq  B^{-1}(sr)\leq k_1B^{-1}(s)B^{-1}(r).\]
\end{prop}
\noindent The above condition, due to \cite[Chapter~II.2.3, Proposition~12]{rao-ren}, is equivalent to the fact that the growth of $B$ is comparable to a certain  power-type function with a fixed exponent.

\medskip

\textit{\underline{Traces and extensions}. } Let us only mention that the theory of the trace operator is not yet extensively examined. There are  results within the doubling setting~\cite{Lacr,Palm} developed in~\cite{Ci-trace} and  further in~\cite{AKKM} towards Orlicz-Slobodetskii spaces.

\subsection{General Musielak-Orlicz setting}\label{ssec:genMOs}

The general isotropic approach is investigated starting from the pioneering monograph of~Nakano~\cite{Nakano} and articles by Skaff~\cite{Sk1,Sk2} and Hudzik~\cite{HH3,HH4}. The preeminent role for the functional analysis of Musielak-Orlicz spaces is played by the  monograph of~Musielak~\cite{Musielak}.
Applications in modelling involving the setting start from classical 
Ball's paper~\cite{ball} on~elasticity, investigated recently e.g. in~\cite{ch-o,Filip}. We refer also to~\cite{gwiazda-non-newt,gwiazda-tmna,gwiazda2,pgasgaw-stokes,Aneta,AWK-unsteady2013} for some developments arising around the theory  of~non-Newtonian fluids and for existence to~some parabolic problems within the setting to~\cite{ASGcoll,ASGpara}. Nowadays intensively investigated fields are also potential theory~\cite{HPHP}, harmonic analysis~\cite{Die-max,HP}, regularity theory~\cite{HPHPAK,hht}, and homogenization within the setting is studied in~\cite{Homog1,Homog2}.  We want to stress embedding results of~\cite{CUH,MaSaSh}. Coming back to PDEs we note that~\cite{ASGpara} and~\cite{pgisazg3} investigate existence of parabolic problems, when modular function depends not only on the space variable, but also on time, that is -- in~the space changing with time. Renormalized solutions to~$L^1$-data problems in nonreflexive anisotropic Musielak-Orlicz spaces are considered in the elliptic setting in~\cite{pgisazg1,gwiazda-ren-ell,gwiazda-ren-cor} and in the parabolic one in~\cite{pgisazg2,pgisazg3,gwiazda-ren-para}. More on this issue can be found in Section~\ref{sec:below}.

\subsubsection*{Basic definitions}

In general, we refrain from precise formulations here, but the following ones have a~fundamental meaning and cannot be omitted. Note that they are all anisotropic.

\begin{defi}\label{def:Nf} Suppose $\Omega\subset\rn$ is an open bounded set. A~function   $M:\Omega\times\rn\to\r$ is called an $N$-function if it satisfies the
following conditions:
\begin{enumerate}
\item $ M$ is a Carath\'eodory function (i.e. measurable with respect to $x$ and continuous with respect to the last variable), such that\begin{itemize}
\item[i)] $M(x,0) = 0$
\item[ii)]  $\inf_{x\in\overline{\Omega}}M(x,\xi)>0\ $ for $\ \xi\neq 0$,
\item[iii)] $M(x,\xi) = M(x, -\xi)$ a.e. in $\Omega$;
\end{itemize}
\item $M(x,\xi)$ is a strictly convex function with respect to $\xi$,
\item $\lim_{|\xi|\to 0}\mathrm{ess\,sup}_{x\in\overline{\Omega}}\frac{M(x,\xi)}{|\xi|}=0$,
\item $\lim_{|\xi|\to \infty}\mathrm{ess\,inf}_{x\in\overline{\Omega}}\frac{M(x,\xi)}{|\xi|}=\infty$.
\end{enumerate}
\end{defi}

\begin{rem} We call an $N$-function $M$ {\it locally integrable} if for any fixed $c\in\rn$ and any measurable set $K\subset\Omega$ we have
\[\int_K M(x,c)dx<\infty.
\]\end{rem}

\begin{defi}\label{def:MOsp} Let $M$ be a locally integrable $N$-function. We deal with three Orlicz-Musielak classes of functions.\begin{itemize}
\item[i)]${\cal L}_M(\Omega)$  --- a generalised Orlicz-Musielak \underline{class} is the set of all measurable functions  $\xi:\Omega\to\rn$ such that
\[\int_\Omega M(x,\xi(x))\,dx<\infty.\]
\item[ii)]${L}_M(\Omega)$  --- a generalised Orlicz-Musielak space is the smallest linear space containing ${\cal L}_M(\Omega)$, equipped with the Luxemburg norm 
\[
||\xi||_{L_M(\Omega)}=\inf\left\{\lambda>0:\int_\Omega M\left(x,\frac{\xi(x)}{\lambda}\right)\,dx\leq 1\right\}.
\]
\item[iii)] ${E}_M(\Omega)$  --- the closure in $L_M$-norm of the set of bounded functions.
\end{itemize}
\end{defi}
\begin{rem} Let us point out  that outside $\Delta_2$-family ${\cal L}_M$ is not a linear space and that space ${E}_M(\Omega)$ is separable. Directly from the definition we   know that
\[{E}_M(\Omega)\subset {\cal L}_M(\Omega)\subset { L}_M(\Omega).\] 
\end{rem}
\begin{defi}
The Young~conjugate $M^\ast$ to a function  $M:\Omega\times\rn\to\r$ is defined by
\[M^\ast(x,\eta)=\sup_{\xi\in\rn}(\xi\cdot\eta-M(x,\xi)),\qquad \eta\in\rn,\ x\in\Omega.\]\end{defi}

\begin{fact}[\cite{Aneta}] If $M$ is a locally integrable $N$-function, then we have the duality \[(E_M(\Omega))^\ast=L_{M^\ast}(\Omega).\]\end{fact}

 \begin{defi}\label{def:convmod}
We say that a sequence $\{\xi_\delta\}_{\delta}$ converges modularly to $\xi$ in~$L_M(\Omega)$ (and denote it by $\xi_\delta\xrightarrow[\delta\to 0]{M}\xi$), if there exists $\lambda>0$ such that
\[
\int_{\Omega}M\left(x,\frac{\xi_\delta-\xi}{\lambda}\right)dx\to 0.
\] \end{defi}

\subsubsection*{Main challenges}
Let us summarise difficulties  resulting from inhomogeneity and Orlicz--type growth obviously applying here. 

\medskip

\textit{\underline{Structure}. }  The considerations on structure of growth and coercivity conditions problems modelled upon 
\begin{equation}
\label{eq:MO}
- \dv A(x,\nabla u)= f \qquad\text{or}\qquad
\pa_t u- \dv A(t,x,\nabla u)= f,
\end{equation}
with $A$ governed by $M$ are basically the same as in the Orlicz setting. In the anisotropic Musielak-Orlicz case under no growth restrictions should read
\begin{equation}
\label{cond:MO:1}A(x,\xi)\xi\geq M(x, \xi )\qquad\text{and}\qquad M^\ast\big(x,A(x,\xi)\big)\leq c M(x, \xi )
\end{equation}
or be merged in one assumption
\begin{equation}
\label{cond:MO:2}c\Big( M(x, \xi)+M^\ast\big(x,A(x,\xi)\big)\Big)\leq A(x,\xi)\xi.
\end{equation}

Let us imprecisely say that $\nabla u$ is considered in $L_M$ and $A(x,\nabla u)$ in $L_{M^\ast}$, where in general $L_{M^\ast}\neq (L_{M})^\ast $. We explain further consequences.

\medskip

Elliptic PDEs can be posed in isotropic Musielak--Orlicz--Sobolev spaces, namely
\[W_0^{m}L_M(\Omega)=\big\{u\in W_0^{1,1}(\Omega):\ \ u\in L_M(\Omega),\ D^\alpha u\in L_M(\Omega),\ |\alpha|\leq m\big\},\]
but in the anisotropic setting there is no generalisation of symmetrization techniques and related embeddings. Thus, as in e.g.~\cite{pgisazg1,gwiazda-ren-ell}, the problems are more likely to be posed in Musielak--Orlicz  spaces having the following structure
\begin{equation}
\label{VM0}
V^M_0(\Omega) =\{\vp\in W_0^{1,1}(\Omega):\ \nabla \vp\in L_M(\Omega)\}.
\end{equation}
Sometimes the meaning of zero boundary value is defined in some other way, e.g. simply by extension by zero outside $\Omega$, see~\cite{CiMa}, or by the closure of smooth compactly supported in $\Omega$ in proper topology, see~\cite{yags}.

Due to the same reason anisotropic parabolic problems would be then considered, as in e.g.~\cite{pgisazg2,pgisazg3,gwiazda-ren-para}, in 
\begin{equation}
\label{VMT0} \begin{split} V_T^M(\Omega)  &=\{u\in L^1(0,T;W_0^{1,1}(\Omega)):\ \nabla u\in L_M(\Omega_T)\}.\end{split}\end{equation}
 
Note that whenever $M$ is a locally integrable $N$-function, each of the spaces $W_0^mL_M(\Omega)$, $V_0^M(\Omega)$, $V_T^M(\Omega)$ is a Banach space, cf.~\cite{Aneta}.

\medskip

\textit{\underline{Lavrentiev's phenomenon}. }  As much as the modular convergence is a~natural tool in the Orlicz setting~\cite{Gossez}, it is still of great significance in Musielak-Orlicz spaces. Here as well weak derivatives are strong ones with respect to the modular topology, but only in absence of~Lavrentiev's phenomenon. 

It is well known in the inhomogeneous setting of~variable exponent spaces, as well as of the double-phase space, that one is equipped with the density of~the~smooth functions only if the interplay between the~behaviour of the modular function with respect to each of~the variables is balanced. This delicate interplay acts on a global level. Indeed, as initially noticed in the papers of Zhikov~\cite{zhikov9798,Zhikov2011}, in~the~variable exponent case the conditions on the exponent $p(x)$ ensuring the absence of~Lavrentiev's phenomenon and the density of smooth functions via mollification, are the same as those required for the regularity of~minimizers. The same phenomenon extends to the several other energies including double phase, see for instance~\cite{yags,bcm-st,ELM,zhikov9798,Zhikov2011}. See Section~\ref{sec:appr} for more details in the general setting.

\medskip

\textit{\underline{Lack of the~integration-by-parts formula}. }  The lack of~the~integration-by-parts formula as in Orlicz spaces, see Section~\ref{sec:orlicz}, is a meaningful difficulty. It is necessary in passing from distributional formulation of~an equation to the particular class of test functions involving the solution itself. Note that in fighting with this the modular approximation plays key role, cf.~\cite{pgisazg2,pgisazg3,gwiazda-ren-para}.

\subsubsection*{Growth restrictions}

We already pointed out how important and structural is the role of the $\Delta_2$--condition in~the Orlicz setting, which is obviously forwarded to Musielak-Orlicz spaces. Let us present a~generalized version of the condition and discuss its role.

The Orlicz definition of $\Delta_2$-condition can be slightly generalized by additive $x$-dependent preturbation. We say that an $N$-function $M:\Omega\times\rn\to\r$ satisfies $\Delta_2$--condition, if there exists a~nonnegative integrable function $h:\Omega\to\r$ such that for some constant $c_{\Delta_2}>0$
\[ M(x,2\xi)\leq c_{\Delta_2}M(x,\xi)+h(x)\qquad\text{for a.e.  } x\in\Omega.
\] 
Note that, as in Orlicz spaces, the family of functions satisfying $\Delta_2$ condition is invariant with respect to multiplications and compositions.

\medskip

\textit{\underline{Doubling examples}. } Let us present examples of non-homogeneous and possibly anisotropic modular functions, for which $M,M^\ast\in\Delta_2$:
\begin{itemize}
\item $M( x,|\xi|)= |\xi |^{p ( x)}$, where $1<< p<<\infty$;  covering variable exponent case with possibly not regular exponent;
\item $M( x,|\xi|)= |\xi |^{p ( x)}\log^{\alpha(x)}(e+|\xi|)$, where $1<< p<<\infty$ and $\alpha\geq0$, or $1\leq p<<\infty$ and $\alpha>>0$;
\item $M( x,\xi)=\sum_i a_i(x)|\xi_i|^{p_i( x)}$, where $1<< p_i<<\infty$, weight functions $ a_i$ are nonnegative and bounded a.e. in $\Omega$, and there is no subset of $\Omega$, where all $a_i$ disappear; this case covers anisotropic weighted variable exponent case with possibly not regular exponent;
\item $M_1(x,|\xi|)=|\xi|^p+a(x)|\xi|^q$ or $M_2(x,|\xi|)=|\xi|^p+a(x)|\xi|^p\log(e+|\xi|)$, where $1< p<q<\infty$ and a weight function $a:\Omega\to [0,\infty)$ is bounded and possibly touching zero; covering the  case of the double-phase space;
\item $M_1(x,|\xi|)=|\xi|^{p(x)}+a(x)|\xi|^{q(x)}$ or $M_2(x,|\xi|)=|\xi|^{p(x)}+a(x)|\xi|^{p(x)}\log(e+|\xi|)$, where $1<< p<q<<\infty$ and a weight function $a:\Omega\to [0,\infty)$ is bounded and possibly touching zero; covering the case of variable exponent double-phase space;
\item $M( x,\xi)=M_0(\xi)+\sum_{i=1}^k  a_i( x)M_i(\xi),$ $k\in \N,$ or  $M( x,\xi)=M_0(\xi)+\sum_{i=1}^N  a_i( x)M_i(\xi_i)$, where the Orlicz modular functions $M_i, M_i^\ast\in\Delta_2$, while the weight functions $ a_i:\Omega\to [0,\infty)$ are bounded and possibly touching zero; covering anisotropic weighted Orlicz case under the most common nonstandard growth conditions.
\end{itemize} 
 
\medskip

\textit{\underline{Non--doubling examples}. } There is a vast range of $N$-functions not satisfying $\Delta_2$-condition, for instance:
\begin{itemize}
\item  $M( x,\xi)=a( x)\left( \exp(|\xi|)-1+|\xi|\right)$  with a bounded   weight $a:\Omega\to (0,\infty)$;
\item  $M( x,\xi)=a( x)|\xi|\log(e+|\xi|)+b(x)|\xi|^p$ with $1<p<\infty$ with bounded and possibly touching zero weights $a,b:\Omega\to [0,\infty)$, if only there is not subset of $\Omega$ of positive measure, where both of them disappear; 
\item $M(x,\xi)=M_1(\xi)+a(x)M_2(\xi)$ with bounded and possibly touching zero weight $a:\Omega\to [0,\infty)$ relating to the double phase space, but with $M_i\not\in\Delta_2$ for $i=1$ or $i=2$. Recall that  $M\not\in\Delta_2$ can be trapped between two power-type functions, see~\cite{CGZG} for a construction;
\item  $M( x,\xi)=a( x)\big(\exp(|\xi_1|)-1\big)+|\xi|^{p(x)}$, $\xi=(\xi_1,\dots,\xi_N)$ with a bounded  and possibly touching zero weight $a:\Omega\to [0,\infty)$ and variable exponent $1<<p<<\infty$ on $\{x: a(x)=0\}$. It is an example of a fully anisotropic modular function (with \underline{not} separated roles of coordinates of the second variable); 
\item  $M( x,\xi)= a( x)|\xi_1|^{p_1( x)}\left(1+|\log|\xi||\right)+\exp(|\xi_2|^{p_2( x)})-1$, when $(\xi_1,\xi_2)\in\R^2$ and $p_i:\Omega\to[1,\infty]$.  It is also an example of a fully anisotropic modular function.
\end{itemize} 

 \medskip

\textit{\underline{The role of growth control}. } The remarkable consequences of the doubling condition, which we want to expose in the beginning are the following facts:\[M\in\Delta_2\quad\implies\quad {E}_M(\Omega)= {\cal L}_M(\Omega)= {L}_M(\Omega),\]
and
 \[E_{M^\ast}  \xlongequal[]{M^\ast\in\Delta_2}L_{M^\ast}\xlongequal[ ]{M \text{ is loc. int.}}(E_{M})^\ast\xlongequal[]{M \in\Delta_2} (L_M)^\ast.\]
Thus, the space equipped with $M,M^\ast\in\Delta_2$ is reflexive and separable, which essentially simplifies methods of PDEs. Then, also the modular topology coincides with the norm one. What is more, if $M$ is  locally integrable, the set of simple functions integrable on $\Omega$ is dense in $L_M(\Omega)$ with respect to~the modular topology,~\cite{Musielak}. Unlike variable exponent Lebesgue spaces or the double phase space generalised Musielak-Orlicz spaces   with $M\not\in\Delta_2$ or $M^\ast\not\in\Delta_2$ stop to be separable or reflexive.

The further consequences of growth control or its lack in the theory of existence of~renormalized solutions to problems of the type~\eqref{eq:MO} are described in Section~\ref{sec:below}.

\section{Problems with data below duality}\label{sec:below}

In studies of partial differential equations of a general form
\[
- \dv A(x,\nabla u)= f \qquad\text{or}\qquad
\pa_t u- \dv A(t,x,\nabla u)= f,\]
when $f$ belongs to the natural dual space to the leading part of the operator, there are already classical existence results. We can mention here those in various settings starting from the classical linear case involving Laplacian e.g.~\cite{stam,lions}, the power-type growth involving $p$-Laplacian e.g.~\cite{ladyzenska,ladyzenska-ell}, ending with those posed in the~Orlicz setting with $\Delta_2$ structure~\cite{Ta0,Ta1} or without growth restrictions~\cite{Gossez3,Gossez2}. In Musielak-Orlicz spaces existence of weak solutions to problems with bounded data are provided in~\cite{pgisazg3,GMWK,pgisazg1,ASGcoll,ASGpara}.

The challenge starts, when one wants to~consider less regular data, namely $f$ merely integrable.

\subsection{Necessity of special notion of solutions}\label{ssec:necessity}
Investigating problems 
\begin{equation} \label{ell-para}
- \dv A(x,\nabla u)= f\in L^1(\Omega)\quad\text{and}\quad
\pa_t u- \dv A(t,x,\nabla u)= f\in L^1(\Omega_T)
\end{equation} 
involving $A$ with a growth governed by a function~from the Orlicz class we need to employ a special notion of solutions. The reason is   that data is merely integrable and it does not belong to the natural dual space to the leading part of the operator. To illustrate the issue let us think first about the classical Poisson equation on   $ \rn$, namely
\[\left\{\begin{array}{rcl}
-\Delta u &=&f,\\
u &=&0.
\end{array}\right.\]
Its solution can be expressed by the mean of the Green function $G$ via the formula
\[u(x)=\int_{\rn} f(y)G(x,y)\,dy \stackrel{ {N>2}}{=}c(N)\int_{\rn} f(y)\frac{1}{|x-y|^{N-2}}\,dy,\]
which apparently does \textbf{not} belong to the natural energy space $W^{1,2}(\Omega)$, when $f\in L^1(\Omega)$. To find a~solution one needs to consider a generalised notion of solutions. The easy way would be analysing distributional solutions, but then we cannot ensure uniqueness. The classical example of Serrin~\cite{Serrin-pat} concerns a linear homogeneous equation of the type
$\dv (A(x)Du) = 0$ defined on a ball, with strongly elliptic and bounded, measurable matrix $A(x)$, that has at least two distributional solutions, among which only  one belongs to the natural energy space $W^{1,2}(B)$. Nonetheless, not accidentally the second solution is called pathological. The point is to distinguish solutions having proper interpretation, say physical interpretation, excluding wild ones. The framework which we need should provide unique and interpretable solutions.

\medskip

Therefore, the key property expected from an interesting special notion of~solutions besides existence is uniqueness. The problem with uniqueness appearing in the linear equation is obviously shared by the $p$-harmonic problem $-\Delta_p u =f\in L^1(\Omega),$ as well as its anisotropic, Orlicz, and Musielak-Orlicz generalisations. Indeed, consider~\eqref{ell-para}, where the growth of $A$ -- the leading part of the operator is governed by $M$ via conditions~\eqref{cond:MO:1} or~\eqref{cond:MO:2}. When on the right--hand side data is merely integrable the weak formulations of~\eqref{ell-para}, i.e.
\[\int_\Omega A(x,\nabla u)\nabla\vp\,dx=\int_\Omega f\vp\,dx,\]
\textbf{cannot} be expected to hold for every  $\vp\in V^M_0(\Omega)$ given by~\eqref{VM0}, respectively, in parabolic case, the weak formulation
\[
-\int_{\Omega_T}(u-u_0)\partial_t \vp \,dx\,dt+\int_{\Omega_T}A(t,x,\nabla u)\cdot \nabla\vp \,dx\,dt=\int_{\Omega_T}f\, \vp \,dx\,dt,\]
\textbf{cannot} hold with every $\vp\in V_T^M(\Omega)$ given by~\eqref{VMT0}.

Let us stress once again that under no~growth conditions on $M$ or $M^\ast$ the understanding of the dual pairing  complicates due to\[(L_M)^\ast\neq L_{M^\ast}\] and even for nice data an approximation result is needed from the very beginning. Indeed, we cannot test an equation by the solution itself (or even its truncation) in order to get  a~priori estimates, because $\nabla u$ is considered in $L_M$ and $A(x,\nabla u)$ is supposed to live in $L_{M^\ast}$. In order to get a priori estimates, we need to test the formulation by a sequence of functions admissible in proper pairings and convergent in the modular topology to the solution. See how it works e.g. in~\cite{pgisazg1,pgisazg2}.

\subsection{Various notions of solutions } \label{ssec:variousnotions} There are a few classical notions of solutions introduced in order to consider a datum $f$ not belonging to the dual space.

DiPerna and Lions investigating the~Boltzmann equation in order to deal with this challenge introduced the notion of \textit{renormalized solutions} in~\cite{diperna-lions} with fundamental developments by Boccardo, Giachetti, Diaz, and Murat~\cite{boc-g-d-m} and Murat~\cite{murat}.  Another seminal idea for problems with $L^1$-data comes from Boccardo and Gallou\"et~\cite{bgSOLA,bgSOLA-cpde}, who studied the \textit{solutions obtained as a~limit of approximation},  SOLA for short. Finally, the \textit{entropy solutions} are considered starting from cornerstones laid by Benilan, Boccardo, Gallou\"et, Gariepy, Pierre, and Vazqu\'ez~\cite{bbggpv}, Boccardo, Gallou\"et, and Orsina~\cite{bgo}, and Dall'Aglio~\cite{dall}. Below we present the definitions in the simplest case presenting clearly the main ideas, i.e. in the case of $p$-harmonic elliptic problem with merely integrable data, namely
\begin{equation}
\label{eq:plap}-\Delta_p u=f\in L^1(\Omega).
\end{equation} It should be stressed from the beginning  that the mentioned  distinct  notions of solutions can coincide. See~\cite{KiKuTu} for an elliptic result and~\cite{DrP} for the equivalence between entropy and renormalized solutions to parabolic problems with polynomial growth and~\cite{ZZ,zhang17} for the corresponding results in the variable exponent and the Orlicz settings, respectively. To our best knowledge there are no such results in Musielak-Orlicz spaces.
 
 \medskip

Before we present the definitions of the three notions of solutions, we  introduce the~symmetric truncation  given by\[
T_k(f)(s)=\left\{\begin{array}{ll}f(s), &\ \text{ if }\ |f(s)|\leq k,\\
k\frac{f(s)}{|f(s)|},&\ \text{ if }\ |f(s)|\geq k .
\end{array}\right.  
\]
We note that according to \cite[Lemma~2.1]{bbggpv}, for every $u\in W^{1,1}(\Omega)$, there exists a unique measurable function $Z_u:\Omega\to\rn$ such that\[\nabla (T_t(u))=\chi_{\{|u|<t\}}Z_u\quad\text{a.e. in } \Omega, \text{ for every }{t>0}.
\]
Thus, in the theory $Z_u$ is called the generalized gradient of~$u$ and, abusing the notation, for $u$ in the space of truncations, it is written simply $\nabla u$ instead of~$Z_u$. 

We consider the space of truncations
\[ {\cal T}^{1,p}(\Omega)= 
 \{u\text{ is measurable in }\Omega : \quad T_k(u)\in W^{1,p}(\Omega) \quad \forall_{k>0}\},
\]
where zero trace can be defined in several ways, for instance as in\[ {\cal T}_0^{1,p}(\Omega)= 
 \{u\text{ is measurable in }\Omega : \quad T_k(u)\in W_0^{1,p}(\Omega) \quad \forall_{k>0}\}.
\] 

We are in position to present the three notions.

\begin{defi}\label{def:sola}
We call a function $u\in W^{1,1}_{loc}(\Omega)$ a \textbf{SOLA} to~\eqref{eq:plap},  if problems \[
- \Delta_p  u_k  = f_k \in L^\infty(\Omega)
\] with  $k\in\N$ and
\[\begin{split}&\lim_{k\to\infty}\int_\Omega \vp\,f_k\,dx=\int_\Omega\vp\,f\,dx \quad\text{ for every function }\ \ \vp\in C_c(\Omega)\\
\text{ and }\qquad &\limsup_{k\to\infty} \int_B f_k \,dx \leq \int_B f \,dx \quad\text{
for every ball }\quad B\subset\Omega,\end{split}\]  
have solutions  $\{u_k\}_k\subset W^{1,p}_{loc}
(\Omega)$ such that \[u_k\xrightharpoonup[k\to\infty]{} u\quad\text{ weakly\ \ in\ \  }W^{1,p }_{loc}
(\Omega).\]
\end{defi}

\begin{defi}\label{def:es}
We call a function $u$ an \textbf{entropy solution} to~\eqref{eq:plap}, when \begin{itemize}
\item[(E1)] $u \in {\cal T}_0^{1,p}(\Omega)$.
\item[(E2)]  for every $\vp\in C_c^\infty(\Omega)$ and for every $k>0$ we have \[\int_{\{|u-\vp|<k\}}|\nabla u|^{p-2}\nabla u\cdot(\nabla u-\nabla\vp)\, dx\leq \int_\Omega T_k(u-\vp)f\,dx\quad\forall_{k>0}\]
\end{itemize} 
\end{defi}

\begin{defi}\label{def:rs}
We call a function $u$ a \textbf{renormalized solution} to~\eqref{eq:plap}, when it satisfies the following conditions.\begin{itemize}
\item[(R1)] $u \in {\cal T}_0^{1,p}(\Omega)$.
\item[(R2)] For every $h\in C^1_c(\R)$ and all $\varphi\in W^{1,p}_0\cap L^\infty (\Omega)$ we have
\[\int_\Omega |\nabla u|^{p-2}\nabla u\cdot\nabla(h(u)\varphi)dx=\int_\Omega fh(u)\varphi\,dx.\]
\item[(R3)] $\int_{\{l<|u|<l+1\}}|\nabla u|^{p}\, dx\to 0$ as $l\to\infty$.
\end{itemize} 
\end{defi}

Stressing the role of the generalized gradient $Z_u$ is of particular meaning when we compare renormalized or entropy solutions to SOLA. Indeed, the notion of SOLA   takes into account only $u\in W^{1,1}_{loc}(\Omega)$, which in the case of~equations involving the $p$-Laplace operator and arbitrary measure data requires $p>2-1/N$. For an explanation see~\cite{bbggpv}.

There are also certain other notions also sharing fundamental property of uniqueness for $L^1$-data. Recently in~\cite{CiMa} in the Orlicz setting the notion of approximable solutions has been introduced, somehow merging the ideas of SOLA and entropy solutions, see also~\cite{ACCZG,CGZG}.

 Some of the mentioned results are relevant in the context of measure data problems. 
 Let us refer to~\cite{bgSOLA,bgSOLA-cpde,bgo,DMMOP,DMOP,Kilp-Maly} and further to e.g.~\cite{BMMP,Dong-Fang,MuPo} for elliptic results and~\cite{boc-ors1,BCM} and further to e.g.~\cite{DrP,P,PPP} for parabolic ones. Nonetheless, the uniqueness in the case of arbitrary measure data is a~long-standing open problem. 

We note that for solutions with data below the natural duality we can provide gradient estimates. As for already classical result we  to~\cite{bgSOLA-cpde}, and further to~\cite{min-grad-est,DHM-p}. The global approach to estimates to measure data problems, that works both in the case of very weak solutions and in the case of energy solutions, has been provided via potential estimates. For a full account of estimates for general elliptic equations of $p$-Laplace-type in the scalar case we refer to~\cite{KuMi-guide}, while for the vectorial one to~\cite{KuMi-vec-npt}, whereas the corresponding results for parabolic problems are provided in~\cite{KuMi-wolff}. For gradient estimates to solutions in generalized settings see~\cite{IC-grad,CiFu,CiMa-ARMA14,CiMa,CiMa-ARMA18,Lieb-aniso}.

\subsection{Towards nonstandard growth and inhomogeneity}

Let us concentrate on the nonstandard growth problems.

\medskip

\textit{\underline{Elliptic existence}. } For the existence results for elliptic problems with nonstandard growth and data below duality we refer  to~\cite{renel2,renel3,Dong,fan12,gwiazda-ren-ell,gwiazda-ren-cor,hhk,le-ex,liuzhao15}. Note that~\cite{BWZ,wit-zim} concern isotropic variable exponent spaces, while~\cite{BeKaSa,BeLaSa,le-ex} the anisotropic ones. Isotropic and reflexive Orlicz spaces are employed in~\cite{renel2,renel3,CiMa,CiMa-ARMA18,Dong-Fang}, isotropic and nonreflexive in~\cite{CGZG}, while anisotropic and nonreflexive ones in~\cite{ACCZG}. In~\cite{fan12,hhk,liuzhao15} isotropic, separable and reflexive Musielak-Orlicz spaces are employed, \cite{Dong}~studies separable, but not reflexive Musielak-Orlicz spaces, 
 while  in anisotropic and non-reflexive Musielak-Orlicz spaces in~\cite{pgisazg1,gwiazda-ren-ell,gwiazda-ren-cor}. 
 
 \medskip
 
\textit{\underline{Parabolic existence}. } For results concerning parabolic problems we refer to~\cite{BWZ,pgisazg2,pgisazg3,gwiazda-ren-para,LG,MMR17,ZZ,zhang17}. Among them the variable exponent setting is employed in~\cite{BWZ,LG,ZZ} and non-reflexive Orlicz-Sobolev spaces are studied in~\cite{MMR17,zhang17}. Problems stated in the~anisotropic and non-reflexive Musielak-Orlicz spaces in~\cite{gwiazda-ren-para} under certain growth conditions on a modular function and in~\cite{pgisazg2,pgisazg3} under regularity restrictions only. Note that  in~\cite{pgisazg3} the investigated space is inhomogeneuos not only with respect to $x$, but also to the time variable.

\medskip

\textit{\underline{The framework}. } Let us concentrate on problems stated in the Musielak-Orlicz setting with data below duality of the type~\[
- \dv A(x,\nabla u)= f\in L^1(\Omega)\qquad\text{and}\qquad
\pa_t u- \dv A(t,x,\nabla u)= f\in L^1(\Omega_T)
\]  under growth and coercivity given by the means of an inhomogeneous and anisotropic $N$-function of general growth via conditions having the form~\eqref{cond:MO:1} or~\eqref{cond:MO:2}. Solutions to such problems  would live in spaces of truncations (cf. Definitions~\ref{def:es}, and~\ref{def:rs})
\[\begin{split} {\cal T}V^M(\Omega) &= \{u\text{ is measurable in }\Omega :  \ T_k(u)\in W^{1,1}_0(\Omega),\  \nabla T_k(u)\in L_M(\Omega) \  \forall_{k>0}\},\\
{\cal T}V_T^M (\Omega) &=\{u\text{ is measurable in }\Omega : \  T_k(u)\in L^1(0,T;W_0^{1,1}(\Omega)):\ \nabla T_k(u)\in L_M(\Omega_T) \ \forall_{k>0}\}.\end{split}\]
Let us recall the notion of generalized gradient $Z_u$ which, abusing the notation,  is usually denoted by $\nabla u$ despite in general $u\not\in W^{1,1}_{loc}(\Omega)$. 

 In the case of renormalized solutions the key property expected from~this type of~solutions, which ensures proper interpretation and uniqueness, describes radiation control. Namely, we expect decay of~energy on the level sets of the solution, which here has a~form
\[\int_{\{l<|u|<l+1\}}A(x,\nabla u)\cdot\nabla u\,dx\xrightarrow[l\to\infty]{}0.\] 
 
\medskip

\textit{\underline{The role of $\Delta_2$-condition}. }  We would like to highlight here the role of the technical assumption $M \in\Delta_2$ and the structural assumption $M^\ast\in\Delta_2$ in studies on~\eqref{eq:MO} and its generalisations in nonreflexive and anisotropic Musielak-Orlicz spaces due to the instance of~\cite{pgisazg1,pgisazg2,pgisazg3,gwiazda-ren-ell,gwiazda-ren-para}. 
 The main idea in each of~these papers involves showing existence of~weak solutions to a bounded regularized problem, then passing to a non-regularized problem with bounded data, and finally to~a non-regularized problem with $L^1$-data. This is a multi-stage  construction of~solutions via passing to the limit starting with regular (smooth) functions. An interesting idea of proofs in the parabolic case in~\cite{pgisazg2,pgisazg3,gwiazda-ren-para} is that the notion of renormalized solution is employed to~get weak solution to~truncated and regularized problem.

 Typically assumption $M\in\Delta_2$ can be imposed for reflexivity, facilitation of computations (e.g. factoring a constant out of a modular function), classical embeddings, or the Aubin Lions Lemma (in order to get almost everywhere convergence). To use the Sobolev embedding it may be enough  instead of~$M\in\Delta_2$ to assume essentially less restrictive condition $M(x,\xi)\geq c|\xi|^{1+\ve}$ for $|\xi|>|\xi_0|$ and arbitrarily small $\ve>0$ (as in~\cite{gwiazda-ren-ell,gwiazda-ren-para}). To obtain almost everywhere convergence one can use the comparison principle (as in~\cite{pgisazg2}). 
 
 Let us stress that resigning from $\Delta_2$-condition on the conjugate of a modular function is much more demanding, since it essentially affects the understanding of the dual pairing. This problem is striking.  While testing equation by $T_k (u)$ (in order to get a priori estimates), assumption $M^\ast\in\Delta_2$ fixes structure of a space. Indeed, we expect the duality pairing: \[\nabla T_k (u)\in L_{M}\quad\text{and}\quad A(x,\nabla T_k (u))\in L_{M^\ast} \xlongequal[]{M^\ast\in\Delta_2}E_{M^\ast}\ \ \text{predual to}\ \   L_{M} ,\] 
followed by~\cite{gwiazda-ren-ell,gwiazda-ren-para}. Otherwise, that is when $M^\ast\not\in\Delta_2$ (i.e. very fast, or slow, or having~irregular growth -- cf.~Section~\ref{sec:orlicz}), to pass it by we need to employ some approximate sequence $(T_k (u))_\delta$ of~admissible functions converging in the modular topology. 
However, since a space is inhomogeneous, we expect Lavrentiev's phenomenon and growth restrictions cannot be just skipped -- they can be traded towards regularity ones. Namely, modular function has to satisfy a condition related to the log-H\"older continuity of a variable exponent capturing decent interplay between a behaviour of a modular function with respect to the space and the gradient variable. See Section~\ref{sec:appr} for more information.

 The proofs in~\cite{pgisazg1,pgisazg2,pgisazg3} are formulated under no growth conditions and in the absence of~Lavrentiev's phenomenon. The method keeps a solution in the modular closure of smooth functions, which coincides with the strong closure when $M,M^\ast\in\Delta_2$. Thus, to avoid considering Lavrentiev's phenomenon, it can be assumed instead that $M,M^\ast\in\Delta_2$ and then an approximation via Mazur's Lemma turns to be sufficient in order to get existence of renormalised solutions. In turn, the mentioned studies include (without additional assumptions) in particular\begin{itemize}
\item[{\it i)}] the anisotropic Orlicz case  under no growth restrictions,
\item[{\it ii)}] reflexive  spaces, that is among others: the variable exponent, the weighted Sobolev, and the double phase space.
\end{itemize}

\section{Density and approximation}\label{sec:appr}

One of important features of the non-homogeneous setting that Musielak-Orlicz spaces inherit is problem with density of smooth functions, which is closely related to~absence of the Lavrentiev phenomenon~\cite{LM}. It is necessary to know how to exclude the possibility of existence of functions which cannot be approximated by regular ones.  

The Lavrentiev phenomenon was already mentioned in the context of variable exponent spaces (with $M(x,s)=|s|^{p(x)}$) with not regular $p(\cdot)$, cf.~\cite{ZV}, and in the context of the double-phase spaces (with $M(x,s)=|s|^{p}+a(x)|s|^{q}$) with $p$ and $q$ are far from each other, cf.~\cite{min-double-reg1}. In brief,  if a Musielak-Orlicz space is equipped with a modular function whose behaviour is not sufficiently balanced, there exist functions that cannot be approximated in the relevant topology by regular (smooth) functions. 
 
The above examples involve reflexive spaces, when the modular topology coicides with the strong one. The strong closure of the smooth functions, however, is not a relevant type of handy approximation here. In fact, as in Orlicz spaces (see Section~\ref{sec:orlicz}) in nonreflexive Musielak-Orlicz spaces  the relevant topology is not the norm topology, but the modular one.  The fundamental results by  Gossez~\cite{Gossez} in the Orlicz setting are extended to~the anisotropic Musielak-Orlicz setting in~\cite{pgisazg1} and refined in the isotropic case in~\cite{yags} under  regularity restrictions on the modular function.

Let us stress again that it entails that kind of the Meyers-Serrin theorem, saying that weak derivatives are strong ones with respect to the modular topology, in Musielak-Orlicz spaces holds only in absence of~Lavrentiev's phenomenon.

\subsection*{Assumptions on the modular function}
To ensure the modular density of smooth functions, it is necessary to impose a restriction  balancing a behaviour of a modular function $M$ for big $|\xi|$ and small changes of~the space variable.
  
In the isotropic case in~\cite{pgisazg1}, the authors prove that it suffices to impose on $M$ continuity condition of log-H\"older-type with respect to $x$, namely for each $\xi\in\rn$ and $x,y,$ such that $|x-y|<\frac{1}{2}$ we have\[
 \frac{M(x,\xi)}{M(y,\xi)}\leq\max\left\{ |\xi|^{-\frac{a_1}{\log|x-y|}}, b_1^{-\frac{a_1}{\log|x-y|}}\right\},\ \text{with some}\ a_1>0,\,b_1\geq 1.
\]  Note that this condition 
 for $M(x,\xi)=|\xi|^{p(x)}$ relates to the log-H\"older continuity condition for a variable exponent $p$, namely there exists $a>0$, such that for $x,y$ close enough and  $|\xi|\geq 1$
\[|p(x)-p(y)|\leq \frac{a}{\log\left(\frac{1}{|x-y|}\right)}.\]
Indeed, whenever $|\xi|\geq 1$
\[ \frac{M(x,\xi)}{M(y,\xi)}= \frac{|\xi|^{p(x)}}{|\xi|^{p(y)}}=|\xi|^{p(x)-p(y)}\leq |\xi|^\frac{a}{\log\left(\frac{1}{|x-y|}\right)}=|\xi|^{-\frac{a}{\log {|x-y|} }}.\]

Let us refer to the approaches of~\cite{hhk,hht} and~\cite{mmos:ap,mmos2013}, where the authors deal with the isotropic modular function of the form $M(x,\xi)=|\xi|\phi(x,|\xi|)$. As for the types of~regularity,  in~\cite{mmos:ap,mmos2013} the authors restrict themselves to the case when \[\phi(x,|\xi|)\le c \phi(y,|\xi|)\qquad\text{when}\qquad |\xi|\in [1,|x-y|^{-n}].\] Meanwhile in~\cite{hhk,hht}, the proposed condition yields \[\phi(x, b|\xi|)\le \phi(y,|\xi|)\qquad\text{when}\qquad\phi(y,|\xi|)\in [1, |x-y|^{-n}].\] In the end we point out that it is explained in~\cite{pgisazg1,yags}, why it was necessary to~fix the previously existing proofs of the approximation theorems whose insufficient conditions got somehow propagated.

\subsection*{Gossez's approximation theorems}

 Let us go to the conditions capturing the known optimal cases, i.e. the log-H\"older condition on the exponent in the variable exponent case and the double-phase space with $\alpha$-H\"older continuous weight $a$ within the sharp range of~the involved parameters. Merging approaches of~\cite{pgisazg1,yags} we consider the modular function $M$ be a locally integrable $N$-function.  We present the approximation result when $M$ with arbitrary growth satisfies anisotropic condition {\it (M)} or much easier to interpret isotropic  condition {\it (M$^{\, i}$)} given below.
\begin{itemize}
\item[{\it (M)}] Consider \[\begin{split} 
 {M}_{x,\varepsilon}(s):=&\inf_{y\in B(x,\varepsilon)}  M(y,s),\qquad\ve>0,\end{split}\]
and $( {M}_{x,\varepsilon})^{\ast\ast}=\big(( {M}_{x,\varepsilon})^{\ast}\big)^{\ast}$ standing for its the second conjugate (coinciding with its greatest convex minorant). Assume that there exists a function  $\Theta:\big[0, {1}/{2} ]\times\rp0\to\rp0$ such that $\Theta(\cdot,s)$ and $\Theta(x,\cdot)$ are nondecreasing functions and for all $x,y\in\overline{\Omega}$ with $|x-y|\leq\frac{1}{2}$ and for any constant $c>0$
\[
M(y, s )\leq\Theta(|x-y|, s )
 ( {M}_{x,\varepsilon})^{\ast\ast}(s)\qquad\text{ with } \qquad\limsup_{\ve\rightarrow0^+}
\Theta^i(\ve, c\ve^{-N})<\infty.
\]  
\item[{\it (M$^{\,i}$)}] Assume that there exists a function  $\Theta^i:\big[0, {1}/{2} ]\times\rp0\to\rp0$ such that $\Theta^i(\cdot,s)$ and $\Theta^i(x,\cdot)$ are nondecreasing functions and for all $x,y\in\overline{\Omega}$ with $|x-y|\leq\frac{1}{2}$ and for any constant $c>0$
\[
M(y, s )\leq\Theta^i(|x-y|, s )M(x, s )\qquad\text{ with } \qquad\limsup_{\ve\rightarrow0^+}
\Theta(\ve, c\ve^{-N})<\infty.
\]
\end{itemize}  

The following approximation result can be understood as excluding the Lavrentiev phenomenon in the class of Musielak-Orlicz spaces equipped with modular functions satisfying regularity condition~{\it (M)}. Let us recall that modular convergence is defined in~Definition~\ref{def:convmod}. 
\begin{theo} \label{theo:a} Assume that $\Omega\subset\rn$ is a bounded Lipschitz domain and  a locally integrable $N$-function $M$ satisfies {\it (M)} or {\it (M$^{\,i}$)}. Then for every $u\in V_0^M(\Omega)\cap L^\infty(\Omega)$, there exist $\lambda>0$ and a~sequence of functions $u_\delta\in C^{\infty}_{c}(\Omega)$ such that $\nabla u_\delta\to\nabla u$   modularly in $ L_M(\Omega)$ for $\delta\to 0$.\end{theo}

The above theorem implies modular approximation result in the variable exponent case under the log-H\"older continuity assumption on an exponent. Recall that  variable exponent spaces (with $1<<p<<\infty$) are reflexive and, consequently, modular and strong closure coincides. Therefore, as a result we end with approximation in the norm topology.

\begin{ex}We have the following isotropic examples of pairs $M$ and $\Theta$ satisfying {\it (M)}, and thus are admissible in our results on density of smooth functions,~\cite{yags}. \label{ex:M-varia}
\begin{enumerate} 
\item \textbf{Orlicz. } If $M(x,s)=M(s)$ is independent of $x$, then it satisfies obviously {\it (M)} by~choosing $\Theta(\tau,s)=1$. Anisotropic case included.

\item  \textbf{Variable exponent. } Suppose that $M(x,s)=|s|^{p(x)}$, $1<<p<<\infty$, satisfies {\it (M)} with $\Theta(\tau,s)= \max\{s^{\omega(\tau)},s^{-\omega(\tau)}\}$, where $\omega(\tau)=c/(\log(1/\tau))$ is modulus of~continuity of $p$. Thus, it is ensured when $p$ is log-H\"older continuous.

\item \textbf{Borderline double-phase. } When $M(x,s)=|s|^p+a(x)|s|^p\log(e+|s|)$, condition {\it (M)} is satisfied with $\Theta(\tau,s)= 1+{\omega(\tau)}\log(e+s^{-N})$, where $\omega(\tau)$ is modulus of~continuity of $a$. For this it is enough to deal with log-H\"older continuous $a$. 

\item  \textbf{Musielak-Orlicz. } Let $M(x,s)=\sum_{i=1}^K k_i(x)M_i(s)+M_0(x,s)$, where for all $i=1,\cdots,K$ there exist functions $\Theta_i:\big[0, {1}/{2} ]\to\mathbb{R}^+$   satisfying \[
k_i(x)\leq\Theta_i(|x-y|)k_i(y)\qquad\text{with}\qquad \lim_{\ve\to 0^+} \Theta_i(\ve )<\infty,\] whereas $M_0(x,s)$ satisfies {\it (M)} with $\Theta_0$. Then, we can take
\[\Theta(\tau, s)=\sum_{i=0}^K \Theta_i(\tau,s).\]
Similar example can be provided with $M(x,s)=\sum_{i=1}^N k_i(x)M_i(s_i)+M_0(x,s)$.
\end{enumerate}
\end{ex}

\bigskip

The method of~\cite{pgisazg1,yags} leads to the sharp result when a modular function has at least power-type growth, then we can relax {\it (M$^{\,i}$)} to the following one (the difference in claim is under the limsup).
\begin{itemize}
\item[{\it (M$^{\,i}_p$) }]  Assume that $M$ satisfies
\[
M(x, s )\geq c|s|^p \qquad\text{with }\  p>1\text{ and }\  c>0
\] and there exists a function  $\Theta^i:\big[0, {1}/{2} ]\times\rp0\to\rp0$ such that $\Theta^i(\cdot,s)$ and $\Theta^i(x,\cdot)$ are nondecreasing functions and for all $x,y\in\overline{\Omega}$ with $|x-y|\leq\frac{1}{2}$ and for~any constant $c>0$
\[
M(y, s )\leq\Theta^i(|x-y|, s )M(x, s )\qquad\text{ with } \qquad\limsup_{\ve\rightarrow0^+}
\Theta(\ve, c\ve^{-\frac{N}{p}})<\infty.
\]
\end{itemize} 

\begin{theo}\label{theo:a-p} Assume that $\Omega\subset\rn$ is a bounded Lipschitz domain and  a locally integrable $N$-function $M$ satisfies {\it (M$_p^{\,i}$)}. Then for every $u\in V_0^M(\Omega)\cap W_0^{1,p}(\Omega)$, there exist $\lambda>0$ and a~sequence of functions $u_\delta\in C^{\infty}_{c}(\Omega)$ such that $\nabla u_\delta\to\nabla u$   modularly in $ L_M(\Omega)$ for~$\delta\to 0$.\end{theo}

There is an interesting application concerning the celebrated case of the Musielak-Orlicz space, when a modular function has at least power-type growth, namely the double-phase space. Let us recall  $H(x,s) =|s|^{p}+a(x)|s|^{q}$.

\begin{ex}[Double phase] \label{ex:d-p} Consider  $1< p<q$ and nonnegative $a\in C^{0,\alpha}_{loc}(\Omega)$ with $\alpha\in(0,1]$, then $M=H$ satisfies the {\it (M$_p^{\,i}$)} with  \[\Theta(\tau,|s|)= C_a\tau^\alpha|s|^{q-p}+1\] with proper limit within the sharp range of parameters, namely
\begin{equation}\label{sharp-range} \frac{q}{p}\leq 1+\frac{\alpha}{N}.\end{equation}
\end{ex}

When we denote the associated space \[
W^{1,H}(\Omega)=\Big\{u\in W_0^{1,1}(\Omega): H(\cdot,|\nabla u|)\in L^1(\Omega) \Big\},
\]
we have the following consequence, which is sharp due to~\cite[Theorem~4.1]{min-double-reg1}.
\begin{rem}\label{rem:double-phase}
Suppose $p,q>1$, $\alpha\in(0,1)$, and $a\in C^{0,\alpha}(\Omega)$ satisfying~\eqref{sharp-range}.  Then for any  $u\in W^{1,p}_0(\Omega)\cap W^{1,H}(\Omega)$  there exist a sequence $\{u_\delta\}_\delta\subset C_c^\infty(\Omega)$ converging to $u$: \[u_\delta\xrightarrow[\delta\to 0]{}u\quad\text{in}\quad W^{1,p}(\Omega)\qquad\text{and}\qquad  \nabla u_\delta\xrightarrow[\delta\to 0]{M} \nabla u\quad\text{modularly in }W^{1,H}(\Omega),\] which entails $H(\cdot,\nabla u_\delta)\xrightarrow[\delta\to 0]{}H(\cdot,\nabla u)$ in $L^{1}(\Omega).$
 \end{rem}  
 
 We have yet another example, for which~\cite{yags} gives the approximation within the prescribed below family of exponents and parameters.
\begin{ex}[Variable exponent double phase] Consider  $1<p_-< p<q<<\infty$ and nonnegative $a\in C^{0,\alpha}_{loc}(\Omega)$ with $\alpha\in(0,1]$, then $M(x,s)=s^{p(x)}+a(x)s^{q(x)}$ satisfies  {\it (M$_p^{\,i}$)}  with  \[\Theta(\tau,|s|)= \max\{s^{\omega_p(\tau)},s^{-\omega_p(\tau)}\}+ \max\{s^{\omega_q(\tau)},s^{-\omega_q(\tau)}\}\left( C_a\tau^\alpha|s|^{\sup_{x\in\Omega}\big(
q(x)-p(x)\big)}+1\right)\] which has the proper limit if \[p,q\ \text{ are log-H\"older continuous \ \ \ and }\quad\sup_{x\in\Omega}\big(q(x)-p(x)\big)\leq\alpha p_-/N.\]
This coincides with the sharp range~\eqref{sharp-range} in the constant exponent case, that is we are then again in Example~\ref{ex:d-p}.
\end{ex} 
 
\subsubsection*{Regularity of domain}

Using the already existing proofs it is easy to provide the above results for broader class of $\Omega$, but only to segment or cone property. 
\begin{op}
How far can we resign from the assumption on a regularity of the boundary $\pa\Omega$ in Theorems~\ref{theo:a} and~\ref{theo:a-p}? 
 \end{op} Let us note that known approximation results  not only in anisotropic Musielak-Orlicz spaces cf.~\cite{pgisazg1,yags}, but even in homogeneous and isotropic setting of Orlicz spaces broader class of domains considered are ones satisfying segment property, cf.~\cite{Gossez}. Relaxing requirement on regularity of $\pa\Omega$ in each of the mentioned settings results would attract attention. 

\bibliographystyle{plain}
\bibliography{pocket-guide-arxiv.bib}

\end{document}